\documentclass[12pt]{amsart}
\usepackage{amsmath}

\newcommand{\oo}{\infty}
\newcommand{\OldestR}{$\mathsf{Oldest}_\mathit{RND}$}
\newcommand{\OldestD}{$\mathsf{Oldest}_\mathit{DET}$}
\newcommand{\N}{\mathbb{N}}

\newcommand{\be}{\begin{enumerate}}
\newcommand{\ee}{\end{enumerate}}
\newcommand{\bi}{\begin{itemize}}
\newcommand{\ei}{\end{itemize}}
\renewcommand{\i}{\item}

\newtheorem{thm}{Theorem}[section]
\newtheorem{prop}[thm]{Proposition}
\newtheorem{prob}[thm]{Problem}

\newtheorem{cor}[thm]{Corollary}
\newtheorem{conj}[thm]{Conjecture}

\theoremstyle{definition}
\newtheorem{rest}[thm]{Restriction}

\title[Random strategies for the Robin Hood game]
{Random strategies with historical memory for the Robin Hood game}
\author{Boaz Tsaban}
\thanks{Supported by the Koshland Center for Basic Research.}
\address{Department of Mathematics,
Weizmann Institute of Science, Rehovot 76100, Israel}
\email{boaz.tsaban@weizmann.ac.il}
\urladdr{http://www.cs.biu.ac.il/\~{}tsaban}

\begin{document}
\begin{abstract}
The \emph{Robin Hood} game is played as follows:
On day $i$, the Sheriff puts $s(i)$
bags of gold in the cave.
On night $i$, Robin removes $r(i)$ bags from the cave.
The game is played for each $i\in\N$. Robin wins if each bag
which was put in the cave is eventually removed from it; otherwise
the Sheriff wins.

Gasarch, Golub, and Srinivasan studied the Robin Hood game
in the case of random strategies where Robin has no historical memory.
We extend their main result to the case of bounded historical memory, and obtain a
hierarchy of provably distinct games.
\end{abstract}

\maketitle

\section{The Robin Hood game}\label{intro}

The \emph{Robin Hood} game $RH(r,s,A)$ is defined for functions
$r,s:\N \to\N $ such that $1\le r(i)<s(i)$ for each $i$, and for a set $A$,
as follows:
\begin{enumerate}
\item On day $i$, the \emph{Sheriff} (of Nottingham) puts $s(i)$ bags
of gold in the cave, each labelled by an element of $A$. No label is used twice
(over the course of the entire game).
\item On night $i$, \emph{Robin} (Hood) removes $r(i)$ bags from the cave.
\end{enumerate}
The game is played for each $i\in\N $.
Robin wins if each bag which was put in the cave is eventually
removed from it; otherwise the Sheriff wins.

It is easy to see that if Robin has an unlimited historical memory
(knowing at each night $i$ which of the bags in the cave appeared
first), then he has a winning strategy: On night $i$ pick $r(i)$
bags out of those which arrived first.

Deterministic strategies for this game were studied,
from the set-theoretic point of view,
in \cite{Sc1, Sc2}.
Gasarch, Golub and Srinivasan \cite{GGS} consider the case
where Robin has \emph{no historical memory}, that is,
he cannot distinguish between the days where the bags were put in the cave.
They suggest the following \emph{probabilistic} strategy for
Robin: On night $i$, remove random $r(i)$ bags out of the cave
(with uniform distribution). They say that Robin wins
\emph{almost surely} if for each bag put in the cave,
its probability of being eventually removed is $1$.
The probability is taken over Robin's coin tosses.
More precisely, the probability that a bag is not removed
is the product of all
probabilities $p_i$ that the bag is not removed at night $i$,
and the Sheriff's winning (or Robin's loosing)
probability is the supremum of all probabilities
$p_x$ that a bag $x$ put in the cave is never removed.\footnote{Our
definition of the Sheriff's winning probability is
simpler than the one given in \cite{GGS},
but both our and the proofs of \cite{GGS} work for both definitions and actually imply that
the definitions are equivalent.}
Let $L(i)=\sum_{j=1}^i s(j)-r(j)$ denote the number of bags in the
cave after night $i$.
The main result in \cite{GGS}
is that Robin wins almost surely if, and only if,
the series
$$\sum_{i=1}^\infty \frac{r(i)}{L(i)+r(i)}$$
diverges; otherwise Robin loses almost surely.

\section{Strategies with bounded historical memory}\label{main}

We generalize the above mentioned result to the case
of bounded historical memory. The typical case is that
Robin can, on each night, identify the
bags put in the cave on the last $k$ days,
where $k$ is constant.
It will turn out that already the natural strategy
for historical memory $k=1$ is strictly stronger than the natural
strategy in the memoryless case. Moreover $k=2$ yields
a strictly stronger strategy than $k=1$, etc.\ (Theorem \ref{domi}).
In fact, the game can be analyzed in a much broader family of cases,
as will be shown in the sequel.

The most general case is that Robin can, on night $i$, identify the
bags put in the cave on the last $b(i)$ days, where $b:\N\to\N\cup\{0\}$
is a function with $b(i)\le i$ for all $i$.\footnote{$b(i)$
stands for the \emph{bound} on Robin's
historical memory on night $i$.}
(So that $b(i)\equiv 0$ is the memoryless
case studied in \cite{GGS}).
It is natural to denote this game by $RH(r,s,b,A)$, but
our analysis below is independent of the set $A$,
so we will simply write $RH(r,s,b)$.
A key observation is that the following natural
restriction leads to a substantial simplification of
the analysis of the generalized games.

\begin{rest}\label{rest}
We pose the restriction that Robin cannot remember anything
that he forgot earlier, that is, $b(i+1)\le b(i)+1$ for each $i$;
equivalently, the function $i-b(i)$ is nondecreasing.
\end{rest}

We suggest the following deterministic and random strategies for Robin,
motivated by the strategy given in \cite{GGS} for the memoryless case:
Call a bag \emph{very old} if Robin cannot tell the day it was put in the
cave. An important observation is that Robin can identify the very old
bags since he can identify the bags which are \emph{not} very old.
\bi
\i[\OldestD:]
On night $i$ Robin chooses any $r(i)$ many bags
out of the very old bags.
If there are less than $r(i)$ many very old bags, then
Robin also chooses some of the oldest bags among the ones he remembers,
so as to choose $r(i)$ bags in total.\footnote{%
To put this more precisely, on night $i$ Robin has a partition of the bags in the cave
into disjoint (possibly empty) sets $S_0,S_1,\dots,S_{r(i)}$ such that $S_0$ is the set of very old bags,
and for $k=1,\dots,r(i)$, $S_k$ is the set of bags put in the cave $r(i)-k$ days ago.
If $|S_0|\ge r(i)$, then Robin chooses any $r(i)$ many bags
out of the bags in $S_0$. Otherwise,
let $m$ be the minimal such that $r(i)<\sum_{k=0}^m|S_k|$.
Then Robin takes all bags in the sets $S_0,\dots,S_{m-1}$, as well as $r(i)-\sum_{k=0}^{m-1}|S_k|$
many bags from $S_m$.
}
\i[\OldestR:] Same as \OldestD{}, but the $r(i)$ bags are chosen
at random, with uniform probability, out of the older bags.\footnote{%
Using the notation of the previous footnote,
If $|S_0|\ge r(i)$, then Robin chooses at random (with uniform probablity) $r(i)$ many bags
out of the bags in $S_0$. Otherwise Robin takes all bags in the sets $S_0,\dots,S_{m-1}$,
as well as $r(i)-\sum_{k=0}^{m-1}|S_k|$ many random bags from $S_m$.
}
\ei
Observe that if \OldestD{} is a winning strategy for Robin, then so is
\OldestR{}.

Write
$$\tilde L(i) = \max\left\{0,\ \sum_{j=1}^{i-b(i)}s(j)-\sum_{j=1}^{i-1}r(j)\right\}$$
Then $\tilde L(i)$ is the number of bags put in the cave on days
$1,2,\dots,i-b(i)$ and not removed until day $i$.

\begin{prop}\label{det}
~\be
\i If $i-b(i)$ is bounded, then \OldestD{} is a winning strategy for Robin
in $RH(r,s,b)$ (for each $r$, and $s$).
\i If there exist infinitely many $i$ such that $\tilde L(i)\le r(i)$, then
\OldestD{} is a winning strategy in $RH(r,s,b)$.
\ee
\end{prop}
\begin{proof}
(1) is easy. To prove (2), assume that a bag was put in the cave on day $d$.
By (1) we may assume that $i-b(i)$ is unbounded.
Let $i$ be such that $d<i-b(i)$. Restriction \ref{rest} ensures that
this will also hold for all larger $i$'s, so we may assume
further that $\tilde L(i)\le r(i)$. This means that on \emph{night} $i$, all bags
put on days $\le i-b(i)$ (in particular, those put on day $d$) were removed from the cave.
\end{proof}

We may therefore make the following additional restriction.
\begin{rest}\label{rest2}
$\tilde L(i)>r(i)$ for all but finitely many $i$.
\end{rest}

\begin{thm}\label{rand}
Assume that $r$, $s$, and $b$ satisfy Restrictions \ref{rest} and \ref{rest2},
and Robin uses the strategy \OldestR{}.
\be
\i If $\sum r(i)/\tilde L(i)=\infty$, then Robin wins almost surely.
\i If $\sum r(i)/\tilde L(i)<\infty$, then the Sheriff wins almost surely.
\ee
\end{thm}
\begin{proof}
If there is $i$ such that $\tilde L(i)\le r(i)$, let $i^*$ be the maximal such $i$.
Then on night $i^*$, all bags put in the first $i^*-b(i^*)$ days are removed,
and $\tilde L(i)>r(i)$ for all $i>i^*$.
Since the convergence of the series in question does not depend on
the first few elements, this shows that we may assume that $\tilde L(i)>r(i)$
for all $i$.

Observe that if $i-b(i)$ is bounded, then $\tilde L(i)$ is eventually equal to $0$,
contradicting our assumption, thus $i-b(i)$ is unbounded.
Assume that a bag was put in the cave on day $d$, then
for each large enough $i$, $d<i-b(i)$
so that on night $i$, the probability that the bag in question is removed
is $r(i)/\tilde L(i)$.

Now, the probability that a bag put in the cave on day $d$
is \emph{not} eventually removed is
\begin{equation}\label{prod}
\prod_{i=d}^\oo\left(1-\frac{r(i)}{\tilde L(i)}\right).
\end{equation}
The product \eqref{prod} converges to $0$ (i.e., Robin wins almost surely)
if $\sum r(i)/\tilde L(i)=\infty$.
If $\sum r(i)/\tilde L(i)<\infty$, then
the product \eqref{prod} is positive for $d=1$.
Thus, its limit when $d\to\infty$ is $1$.
Consequently, the Sheriff wins almost surely.
\end{proof}

Note that in Theorem \ref{rand}, (1) implies that the other direction of (2)
also holds, and (2) implies that the other direction of (1)
also holds. Thus, the theorem gives an exact characterization of
when Robin wins almost surely and when the Sheriff wins almost surely,
and shows that it is always the case that one of them wins almost surely.
In the case $b(i)\equiv 0$, Theorem \ref{rand} reduces
to the main result of \cite{GGS}, described at the end of Section \ref{intro}.

\section{Historical memory helps}

To make sure that the generalization made in \ref{rand} is not trivial,
we must find instances where additional historical memory changes the strategy's status
from a loosing strategy to a winning one.

Assume that $b,c:\N\to\N$. We say that $c$ \emph{eventually dominates}
$b$ if there exists $m$ such that for all $i>m$, $b(i)<c(i)$.
\begin{thm}\label{domi}
Assume that $c$ eventually dominates $b$, $i-c(i)$ is unbounded, and $b$ and $c$
satisfy Restrictions \ref{rest} and \ref{rest2}.
Then there exist functions $r,s:\N\to\N$ such that for Robin,
\OldestR{} is a (surely) winning strategy in $RH(r,s,c)$ and
an almost-surely loosing strategy in $RH(r,s,b)$.
\end{thm}
\begin{proof}
Since convergence of series does not depend on the first few terms,
we may assume that for each $i$, $b(i)<c(i)$.
Since $\tilde L(i)$ also depends on the bounding function $b$,
let us denote it here by $\tilde L_b(i)$.
It follows that for each
$r$ and $s$, $\tilde L_{c}(i) < \tilde L_{b}(i)$ for all $i$.
It thus suffices to consider the case where $c(i)=b(i)+1$ for each $i$,
and therefore $\tilde L_c(i)+s(i-b(i)) = \tilde L_b(i)$.

At step $i$ of the construction we have the definition of
$s$ at $1,\dots, i-c(i)$ and $r$ at $1,\dots, i-1$,
and therefore $\tilde L_c(i)$ is defined.
Define $r(i)=\max\{i,\tilde L_c(i)\}$,
and define $s$ on $i-c(i)+1,\dots,i+1-c(i+1)$ to be $r(i)^3$
(so that $r(i)/\tilde L_{b}(i)=r(i)/(\tilde L_{c}(i)+s(i-b(i))) = r(i)/(r(i)+r(i)^3) <1/r(i)^2\le 1/i^2$).

Since $r(i)\ge\tilde L_c(i)$ for each $i$, we have by Proposition \ref{det} that
\OldestR{} is a winning strategy in $RH(r,s,c)$.
Now, $\sum r(i)/\tilde L_b(i)\le\sum 1/i^2<\infty$, so by
Theorem \ref{rand}, \OldestR{} is an almost-surely loosing strategy in $RH(r,s,b)$.
\end{proof}

In particular, we have the following.

\begin{cor}
For each $n=0,1,2,\dots$,
there exist functions $r,s:\N\to\N$ such that for Robin,
\OldestR{} is a (surely) winning strategy in $RH(r,s,n+1)$ and
an almost-surely loosing strategy in $RH(r,s,n)$.
\end{cor}

\section{Random Sheriff}

The authors of \cite{GGS} pose the question of the behavior of
the Robin Hood game when the Sheriff's strategy is random as well.
Our analysis in Section \ref{main}, being independent of the Sheriff's
moves, shows that the results apply to this case as well.
In addition to the first strategy of \cite{GGS} which was described in
the introduction to the present paper (Section \ref{intro}),
a second random strategy for Robin is sketched in \cite{GGS}, and
is conjectured to be an almost-surely winning strategy
in the game $RH(1,s,[0,1])$ where $s$ is constant.

While we are unable to analyze the second strategy for lack of
some details, we can see that the first strategy already
works.
Here $b(i)=0$ (no historical memory) and
$L(i)=\sum_{j=1}^i s(j)-r(j)=\sum_{j=1}^i s-1= i(s-1)$,
therefore
$$\sum_{i=1}^\oo \frac{r(i)}{\tilde L(i)} = \sum_{i=1}^\oo \frac{r(i)}{L(i)+r(i)} =
\sum_{i=1}^\oo \frac{1}{L(i)+1} = \sum_{i=1}^\oo \frac{1}{i(s-1)} = \oo,$$
thus by Theorem \ref{rand}, \OldestR{} is an almost-surely winning
strategy in this game.
Since \OldestR{} coincides with the first strategy
of \cite{GGS}, we have that the first strategy of \cite{GGS} is also
an almost-surely winning strategy against a \emph{random} Sheriff in
$RH(1,s,[0,1])$ (as well as any game $RH(r,s,A)$ with
$\sum r(i)/(L(i)+r(i))$ diverging).
Theorem \ref{rand} is an extension of this phenomenon to the case of
nonzero historical memory.

\section{Open problems}

Among the problems which naturally arise, the following two seem
to be the most interesting.

\begin{conj}
\OldestR{} is the best random strategy among the strategies which are independent of
the index set $A$.
\end{conj}

Our analysis repeatedly uses Restriction \ref{rest} posed on the
bounding function $b$.
\begin{prob}
Analyze the the case where $b$ does not satisfy Restriction \ref{rest}.
\end{prob}

Finally, our strategy \OldestR{} uses unboundedly many ``random bits''.
The referee suggests the problem whether good strategies exist,
in which the number of random bits used at each specific step is bounded
by some constant.


\begin{thebibliography}{00}
\bibitem{GGS}
W.\ Gasarch, E.\ Golub, and A.\ Srinivasan,
\emph{When does a random Robin Hood win?},
Theoretical Computer Science \textbf{304} (2003),
477--484.

\bibitem{Sc1}
M.\ Scheepers,
\emph{Variations of a game of Gale (ii): Markov strategies},
Theoretical Computer Science \textbf{129} (1994),
385--396.

\bibitem{Sc2}
M.\ Scheepers and W.\ Weiss,
\emph{Variations of a game of Gale (iii): remainder strategies},
Journal of Symbolic Logic \textbf{62} (1997),
1253--1264.

\end{thebibliography}
\end{document}